\documentclass[12pt,a4paper]{article}
\usepackage[top=1in, bottom=1in, left=1in, right=1in]{geometry}
\usepackage[utf8]{inputenc}
\usepackage{amsmath}
\usepackage{textcomp}
\usepackage{comment}
\usepackage{amssymb}
\usepackage{diagbox}
\usepackage{ytableau}
\usepackage{wrapfig}

\usepackage{youngtab}
\usepackage{graphicx,epstopdf}
\usepackage{ mathrsfs }
\usepackage{graphics}
\usepackage{stmaryrd}
\usepackage{amsfonts}
\usepackage{hyperref}
\usepackage{amsthm}
\usepackage[skip=2pt]{caption}
\usepackage[figurename=Fig.]{caption}
 \usepackage[usenames,dvipsnames]{pstricks}
 \usepackage{epsfig}
 \usepackage{pst-grad} 
 \usepackage{pst-plot} 
 \usepackage{algorithm}
 \usepackage[noblocks]{authblk}
\usepackage{wasysym}
\usepackage[utf8]{inputenc}

\newtheorem{thm}{Theorem}[section]
\newtheorem{thm*}{Theorem}[section]

\newtheorem{lem*}{Lemma}[section]
\theoremstyle{definition}

\title{Circulant graphs}

\begin{document}
\title{\textbf{Total Colorings of Some Classes of Four Regular Circulant   Graphs}}

\author[1]{\ R. Navaneeth}

\author[1]{\ J. Geetha}

\author[1]{\ K. Somasundaram}

\author[2]{Hung-Lin Fu}

\affil[1]{Department of Mathematics, Amrita School of Engineering-Coimbatore\\ Amrita Vishwa Vidyapeetham, India.

\{r\_navaneeth, j\_geetha, s\_ sundaram\}@cb.amrita.edu}
\affil[2]{Department of Applied Mathematics, National Yang Ming Chiao Tung University, Hsinchu 30010, Taiwan. 

hlfu@math.nctu.edu.tw}

\date{}

\maketitle
\begin{abstract}
\indent	 The total chromatic number, $\chi''(G)$ is the minimum number of colors which need to be assigned to obtain a total coloring of the graph $G$. The Total Coloring Conjecture (TCC) made independently by Behzad and Vizing  that for any graph, $\chi''(G) \leq \Delta(G)+2 $, where $\Delta(G)$ represents the maximum degree of $G$. In this paper we obtained the total chromatic number for some classes of  four regular circulant graphs.
\end{abstract}

	\section{Introduction}
	
	Let $G$ be a simple graph with vertex set $V(G)$ and edge set $E(G)$. The \textit{total coloring} of a graph $G$ is an assignment of colors to vertices and edges such that no two adjacent vertices or edges or edges incident to a vertex receives a same color. The \textit{total chromatic number} of a graph $G$, denoted by $\chi''(G)$, is the minimum number of colors required for its total coloring. It is clear that $ \chi''(G) \geq \Delta (G)+1 $, where $\Delta(G)$ is the maximum degree of \textit{G}.\\
	
	Behzad \cite{BEZ} and Vizing \cite{VGV} have independently proposed the Total Coloring Conjecture (TCC) which states that any simple graph $G$, $\chi''(G)\leq{\Delta(G)+2}$. The graphs that can be totally colored by at-least $\Delta(G)+1$ colors are said to be Type I graphs whereas the graphs which can be colored by $\Delta(G)+2$ colors are said to be Type II graphs. The decidability algorithm for total coloring is NP-complete even for cubic bipartite graph \cite{10}.  Good survey of techniques and other results on total coloring can be found in Yap \cite{YAP}, Borodin \cite{OVB} and Geetha et al. \cite{GNS}. In this paper, we obtain the total chromatic number of some four regular circulant graphs are Type I.\\
	
	\section{Four Regular Circulant Graphs}
	
	For a sequence of positive integers $1\le{d_1}<d_2<...<d_l\le{\lfloor\frac{n}{2}\rfloor}$, the circulant graph $G=C_n(d_1,d_2,...,d_l)$ has vertex set $V(G)=\{0,1,2,...,n-1\}$ and two vertices $x$ and $y$ are adjacent if $x\equiv{(y\pm{d_i})\mod n}$ for some $i$ where $1\le{i}\le{l}$. A power of cycles graph $C_n^k$ is a graph with vertex set $V(G)=\{0,1,2,...,n-1\}$  and two vertices $x$ and $y$ are adjacent if and only if $|x-y| \leq k$.  It is easy to see that the four regular circulant graph $C_n(1,2)\cong C_n^2$. Campos and de Mello \cite{6}  proved that $C^2_n, n\neq 7$, are Type I and $C^2_7$ is Type II.  We know that $K_{4,4}$ is a four regular circulant graph and it is Type II \cite{YAP}. A Unitary Cayley graph is a circulant graph with vertex set $V(G)=\{0,1,2,...,n-1\}$ and two vertices $x$ and $y$ are adjacent if and only if $gcd((x-y),n)=1$. Prajnanaswaroopa et al.\cite{PGS}, proved that al most all Unitary Cayley graphs of even order are Type I and odd order satisfies TCC. 
	In this paper, we considered the four regular circulant graphs of the form $C_n(a,b), 1 \leq a < b \leq {\lfloor\frac{n-1}{2}\rfloor}$. \\
	
	Mauro and Diana ~\cite{Mau} proved that the graphs $C_n(2k,3)$ are Type I for $n = (8 \mu  + 6 \lambda )k$, with $k\geq 1$ and non-negative integers $\mu$ and $\lambda$. Riadh Khennoufa and  Olivier Togni~\cite{kt} studied total colorings of circulant graphs and proved that every 4-regular
    circulant  graphs  $C_{6p}(1,k),$ $p\geq 3$ and $k< 3p$ with $k\equiv 1 \mod 3$ or $k\equiv 2 \mod 3$, are Type I. Other cases are still open.\\
    
    Also they proved that the total chromatic number of $C_{5p}(1,k)$, with $k<\frac{5p}{2}, k \equiv 2 \mod 5, k\equiv 3 \mod 5 $ is 5. In the following theorem, we  prove that the graphs $C_{5p}(1,k)$ are Type I for the remaining cases $k\equiv 1 \mod 5$ and  $k\equiv 4 \mod 5$.
    
\begin{thm}
	The circulant graphs $C_{5p}(a,b)$, where $  p\ge 1, a, b \not\equiv {0\mod5}$, are Type I.

\end{thm}
\begin{proof}
	Let $q_1=gcd(5p,a)$ and $q_2=gcd(5p,b)$. The circulant graphs $G=C_{5p}(a,b)$ are four regular graphs with $q_1$ cycles of order $\frac{5p}{q_1}$ and $q_2$ cycles of order $\frac{5p}{q_2}$. \\
	Let $\varphi:V(G)\cup{E(G)\rightarrow{C=\{0,1,2,3,4\}}}$ be a mapping of $G$ defined as follows.\\
	The vertices $v_i$ of $C_{5p}(a,b)$ can be colored by $\varphi(v_i)=i\mod5, 0\le{i}\le{5p-1}$. \\
	Since, the all the four graphs $C_{5p}(a,b)$, $C_{5p}(b,a)$,  $C_{5p}(n-a,b)$ and $C_{5p}(a,n-b)$ are isomorphic to each other,  for the edge colorings, we need to consider the three cases.\\	
	
	\noindent Case 1: $a\equiv{1\mod5}$ and $b\equiv{1\mod5}$.  \\
	The edges of cycles can be colored by setting $\varphi(v_iv_{(i+a)\mod5})=(i+2)\mod5$ and  $\varphi(v_iv_{(i+b)\mod5})=(i+4)\mod5$. If $\varphi(v_i)=c$ where $c\in{C}$, then  $\varphi(v_iv_{(i+a)\mod5})=(c+2)\mod5$, $\varphi(v_{(i-a)\mod5}v_i)=(c+1)\mod5$, $\varphi(v_iv_{(i+b)\mod5)}=(c+4)\mod5$ and $\varphi(v_{(i-b)\mod5}v_{i})=(c+3)\mod5$. \\

	\noindent Case 2: $a\equiv{2\mod5}$ and $b\equiv{2\mod5}$.  \\
	The edges of cycles can be colored by setting $\varphi(v_iv_{(i+a)\mod5})=(i+3)\mod5$ and  $\varphi(v_iv_{(i+b)\mod5})=(i+4)\mod5$. If $\varphi(v_i)=c$ where $c\in{C}$, then  $\varphi(v_iv_{(i+a)\mod5})=(c+3)\mod5$, $\varphi(v_{(i-a)\mod5}v_i)=(c+1)\mod5$, $\varphi(v_iv_{(i+b)\mod5)}=(c+4)\mod5$ and $\varphi(v_{(i-b)\mod5}v_{i})=(c+2)\mod5$.\\

	\noindent Case 3: $a\equiv{1\mod5}$ and $b\equiv{2\mod5}$.  \\
	 The edges of cycles can be colored by setting $\varphi(v_iv_{(i+a)\mod5})=(i+3)\mod5$ and  $\varphi(v_iv_{(i+b)\mod5})=(i+1)\mod5$. If $\varphi(v_i)=c$ where $c\in{C}$, then  $\varphi(v_iv_{(i+a)\mod5})=(c+3)\mod5$, $\varphi(v_{(i-a)\mod5}v_i)=(c+2)\mod5$, $\varphi(v_iv_{(i+b)\mod5)}=(c+1)\mod5$ and $\varphi(v_{(i-b)\mod5}v_{i})=(c+4)\mod5$.\\ Therefore, five colors are used for totally color the graph.
	
	\end{proof}

	In the following theorem, we prove some classes of four regular circulant graphs $C_n(a,b)$ of order $3p$ are Type I.
	
\begin{thm}
	Let $p$ be an odd integer. Then circulant graphs $C_{3p}(a,b)$ with $gcd(a,b)=1$ and $\frac{3p}{gcd(3p, b)}=3s, s \in N$ are Type I.
	
\end{thm}

\begin{proof}

	Let $q_1=gcd(3p,a)$ and $q_2=gcd(3p,b)$. The circulant graphs $G=C_{3p}(a,b)$ are four regular graphs with $q_1$ cycles of order $\frac{3p}{q_1}$ and $q_2$ cycles of order $\frac{3p}{q_2}$.  Let $\varphi:V(G)\cup{E(G)\rightarrow{\{0,1,2,3,4\}}}$ be a mapping obtained by the following process.\\
	
	Let $C_{i}$ be the cycles of order $\frac{3p}{q_2}$ with the vertices $v_{ia}$, $0\le{i}\le{q_{2}-1}$. First we consider the cycle $C_0$. If $q_2=1$ then the vertices and edges of $C_{0}$ are colored with the colors $0,3$ and $1$ cyclically, starting with $v_{0}$ receiving the color $0$. Otherwise the vertices and edges of $C_{0}$ are colored with the colors $3,1$ and $0$ cyclically, starting with $v_{0}$ receiving the color $1$. Now, consider the cycle $C_1$. The vertices and edges of $C_{1}$ are colored with the  colors $1,0$ and $4$ cyclically, starting with $v_{a}$ receiving the color $1$. For the cycles $C_i$,  $2\le{i}\le{q_{2}-1}$, if $i$ is even then the vertices and edges  of $C_{i}$  are colored with the colors $0,2$ and $1$ cyclically, starting with $v_{ia}$ receiving the color $0$ and $i$ is odd they are colored with the colors $1,0$ and $2$ cyclically, starting with $v_{ia}$ receiving the color.\\
	
	The edges of cycles of order $\frac{3p}{q_1}$  are colored in the following way:  if vertex $v_{i}\in{C_{0}}$ then $\varphi(v_{i}v_{i+a})=2$, if $v_{i}\in{C_{i}}$ where $i$ is odd then $\varphi(v_{i}v_{i+a})=3$ and if $v_{i}\in{C_{i}}$ where $i$ is even then $\varphi(v_{i}v_{i+a})=4$. Therefore, only five colors are used for totally coloring the graph. Hence, $\varphi$ is a Type I coloring of $G$
	
\end{proof}

For the circulant graphs $C_n(a,b)$ where $n=3p$ and $p$ is odd, which we considered in the above theorem, the value of $b$ is restricted to factors and multiple of $p$. In the following theorem, we consider few classes of circulant graphs $C_{n}(1,k)$ where $n=9p$ are Type I.	

  \begin{thm}
	For each $k\in{\{2,3,...,{\lfloor\frac{9p-1}{2}\rfloor}\}}$, every circulant graph $C_{9p}(1,k)$ with $\frac{9p}{gcd(9p,k)}=3s, s \in N$ is Type I.
  \end{thm}

  \begin{proof}
Let $q=gcd(9p,k)$. The circulant graphs $G=C_{9p}(1,k)$ are four regular graphs with $q$ internal  cycles of order $\frac{9p}{q}$, which are disjoint, and one outer cycle of order $9p$. Let $\varphi:V(G)\cup{E(G)\rightarrow{\{0,1,2,3,4\}}}$ be a mapping obtained by the following process.\\

  \noindent Case 1: $p$ is even.\\

 When $p$ is even, $9p$ will be a multiple of 6. Riadh Khennoufa and  Olivier Togni~\cite{kt} proved that the total chromatic number of $G=C_{6p}(1,k)$, with $k<\frac{5p}{2}, k \equiv 1 \mod 3, k\equiv 2 \mod 3 $ is 5. 	From this, one can easily see that the circulant graph $C_{9p}(1,k)$, where $p\ge 1$ and  $k<\frac{9p}{2}, k \equiv 1 \mod 3, k\equiv 2 \mod 3 $ are Type I. 

  Now, we consider the remaining case, $k \equiv 0 \mod 3$.\\
      The vertices $v_i$  are colored  by $\varphi(v_i)=(i\mod3+{\lfloor{\frac{i}{k}}\rfloor}\mod3)\mod3$ if $q=k$, else the vertices $v_i$ are colored by $\varphi(v_i)=(2i\mod3-{\lfloor{\frac{i}{3}}\rfloor}\mod3)\mod3$. The edges of the internal cycles can be colored by setting 
      $\varphi(v_{i} v_{(i+k)\mod9p})=(2\varphi(v_{(i+k)\mod9p})-\varphi(v_{i}))\mod3$. In this coloring process,  the vertices and the internal edges of $G$ are colored  using only three colors $0,1$ and $2$. Now, the edges of the outer cycle can be colored with two colors $3$ and $4$ as it is of even order. Therefore, five colors are used for total coloring the graph, hence the graph $C_{9p}(1,k)$ is Type I. \\
	
 \noindent Case 2: $p$ is odd.\\
 \\
 The case when $q\neq{1}$, follows from Theorem 2.2. Now, we consider the case when $q={1}$.
  	\\
  	Sub case 2.1: $k\equiv{1\mod9}$\\
  The vertices $v_i$ where $0\le{i}\le{9p-1}$ of $G$ can be colored by $\varphi(v_i)=i\mod3+(\lfloor{\frac{i}{3}}\rfloor\mod3)\lfloor{\frac{i\mod3}{2}}\rfloor$. The edges of the internal cycles can be colored by setting if $i\equiv{1\mod3}$ then  $\varphi(v_{i}v_{(i+k)\mod9p})=\varphi(v_{(i+k)\mod9p})+1-3\lfloor{\frac{\varphi(v_{(i+k)\mod9p}v_i)\mod5}{5}\rfloor}$ else by $\varphi(v_{i}v_{(i+k)\mod9p})=\varphi(v_{(i+2k)\mod9p})$. The colors used for vertex $v_i$ and edges incident to it is given in the table below as an ordered triplet $(\varphi(v_i),\varphi(v_{(i-k)\mod9p}v_{i}),\varphi(v_{i}v_{(i+k)\mod9p}))$. 
  \begin{center}
  		\begin{tabular}{|c|c|c|c|c|c|c|c|c|}
  			\hline
  			\textbf{$x=0$} &
 			\textbf{$x=1$} & 
  			\textbf{$x=2$} & 
  			\textbf{$x=3$} & 
  			\textbf{$x=4$} & 
  			\textbf{$x=5$} & 
  			\textbf{$x=6$} & 
  			\textbf{$x=7$} & 
  			\textbf{$x=8$} \\				
  			\hline
 			\textbf{$(0,1,2)$} &
  			\textbf{$(1,2,3)$} & 
  			\textbf{$(2,3,1)$} & 
  			\textbf{$(0,1,3)$} & 
  			\textbf{$(1,3,4)$} & 
  			\textbf{$(3,2,1)$} & 
  			\textbf{$(0,1,4)$} & 
  			\textbf{$(1,4,2)$} & 
  			\textbf{$(4,2,1)$} \\
  			\hline
  		\end{tabular}
  		$x\equiv{i\mod9}$
  	\end{center}
  The common missing color for vertices $v_i$ and $v_{i+1}$ can be used for coloring the edge $v_{i}v_{i+1}$. Therefore, five colors are used for totally coloring the graph, hence $G=C_{9p}(1,k)$ Type I.\\
  \\
  Sub case 2.2: $k\equiv{4\mod9}.$\\
  The vertices $v_i$ where $0\le{i}\le{9p-1}$ of $G$ can be colored by $\varphi(v_i)=i\mod3+(\lceil{\frac{i}{3}}\rceil\mod3)\lfloor{\frac{i\mod3}{2}}\rfloor$. The edges of the internal cycles can be colored by setting if $i\equiv{1\mod3}$ then  $\varphi(v_{i}v_{(i+k)\mod9p})=\varphi(v_{(i+k)\mod9p}+1)$ else by $\varphi(v_{i}v_{(i+k)\mod9p})=\varphi(v_{(i+2k)\mod9p})$. The colors used for vertex $v_i$ and edges incident to it is given in the table below as an ordered triplet $(\varphi(v_i),\varphi(v_{(i-k)\mod9p}v_{i}),\varphi(v_{i}v_{(i+k)\mod9p}))$
 	\begin{center}
  		\begin{tabular}{|c|c|c|c|c|c|c|c|c|}
  			\hline
  			\textbf{$x=0$} &
  			\textbf{$x=1$} & 
  			\textbf{$x=2$} & 
  			\textbf{$x=3$} & 
  			\textbf{$x=4$} & 
  			\textbf{$x=5$} & 
  			\textbf{$x=6$} & 
  			\textbf{$x=7$} & 
  			\textbf{$x=8$} \\				
  			\hline
  			\textbf{$(0,1,2)$} &
  			\textbf{$(1,4,2)$} & 
  			\textbf{$(3,4,1)$} & 
  			\textbf{$(0,3,1)$} & 
  			\textbf{$(1,2,3)$} & 
 			\textbf{$(4,2,1)$} & 
 			\textbf{$(0,1,4)$} & 
 			\textbf{$(1,3,4)$} & 
 			\textbf{$(2,3,1)$} \\
  			\hline
  		\end{tabular}
 		$x\equiv{i\mod9}$
  \end{center}
  The common missing color for vertices $v_i$ and $v_{i+1}$ can be used for coloring the edge $v_{i}v_{i+1}$. Therefore, five colors are used for totally coloring the graph, hence $G=C_{9p}(1,k)$ Type I.\\
 \\
   Sub case 2.3: $k\equiv{7\mod9}.$\\
  	The vertices $v_i$ where $0\le{i}\le{9p-1}$ of $G$ can be colored by $\varphi(v_i)=i\mod3+(\lfloor{\frac{i}{3}}\rfloor\mod3)\lfloor{\frac{i\mod3}{2}}\rfloor$. The edges of the internal cycles can be colored by setting $\varphi(v_{i}v_{(i+k)\mod9p})=\varphi(v_{(i+k+1)\mod9p})$. The colors used for vertex $v_i$ and edges incident to it is given in the table below as an ordered triplet $(\varphi(v_i),\varphi(v_{(i-k)\mod9p}v_{i}),\varphi(v_{i}v_{(i+k)\mod9p}))$. 
  \begin{center}
	
  		\begin{tabular}{|c|c|c|c|c|c|c|c|c|}
  			\hline
 			\textbf{$x=0$} &
  			\textbf{$x=1$} & 
 			\textbf{$x=2$} & 
		\textbf{$x=3$} & 
 			\textbf{$x=4$} & 
  			\textbf{$x=5$} & 
  			\textbf{$x=6$} & 
  			\textbf{$x=7$} & 
  			\textbf{$x=8$} \\				
  			\hline
  			\textbf{$(0,1,4)$} &
  			\textbf{$(1,2,0)$} & 
  			\textbf{$(2,0,1)$} & 
  			\textbf{$(0,1,2)$} & 
  			\textbf{$(1,3,0)$} & 
  			\textbf{$(3,0,1)$} & 
  			\textbf{$(0,1,3)$} & 
  			\textbf{$(1,4,0)$} & 
  			\textbf{$(4,0,1)$} \\
  			\hline
  		\end{tabular}
  		$x\equiv{i\mod9}$
  	\end{center}
 The common missing color for vertices $v_i$ and $v_{i+1}$ can be used for coloring the edge $v_{i}v_{i+1}$. Therefore, five colors are used for the total coloring the graph, hence $G=C_{9p}(1,k)$ is Type I.\\
\end{proof}	
	
In Theorem 2.2, we considered few classes of circulant graph $C_n(a,b)$ of order $n=3p$, where $p$ is an odd integer. Now, in the following theorem, we consider few classes of four regular circulant graphs $C_n(a,b)$ of order $n=6p$.
\begin{thm}
Every circulant graph $C_{6p}(a,b)$ where $a,b\not\equiv{0\mod3}$ is Type I, if $p$ is even. Also, $C_{6p}(a,b)$ where $a,b\not\equiv{0\mod3}$ is Type I, if $p$ is odd and $gcd(a,b)=1$.
    \end{thm}
    
 \begin{proof}
    	Let $q_1=gcd(6p,a)$ and $q_2=gcd(6p,b)$. The circulant graphs $G=C_{6p}(a,b)$ are four regular graphs with $q_1$ cycles of order $\frac{6p}{q_1}$ and $q_2$ cycles of order $\frac{6p}{q_2}$. The circulant graphs considered in the hypothesis can be colored in a similar manner irrespective of $p$ being odd or even, if $\frac{6p}{q_1}$ is odd we swap the value of $a$ and $b$, as graph $C_n(a,b)$ is isomorphic to $C_n(b,a)$. Let $\varphi:V(G)\cup{E(G)\rightarrow{\{0,1,2,3,4\}}}$ be a mapping.\\
    The vertices $v_i$  are colored  by $\varphi(v_i)=i\mod3$ and the edges of cycles of order $\frac{6p}{q_2}$ be colored by setting 
      $\varphi(v_{i} v_{(i+a)\mod6p})=(2\varphi(v_{(i+a)\mod3p})-\varphi(v_{i}))\mod3$. Now, the edges of cycle $\frac{6p}{q_1}$ with two colors $3$ and $4$ as it is a cycle with even order. Therefore, five colors are used for the total coloring $\varphi$ of the graph, hence graph $G$ is Type I.
    \end{proof}

\end{document}